\documentclass[11pt, letterpaper,english]{amsart}
\usepackage[left=1in,right=1in,bottom=1in,top=1in]{geometry}
\usepackage{amsmath,amsthm,amssymb}
\usepackage{mathrsfs}
\usepackage[all,cmtip]{xy}
\usepackage{cite}
\usepackage{hyperref}
\usepackage{enumerate}
\usepackage{graphicx}

\theoremstyle{plain}
\newtheorem{lem}{Lemma}

\newtheorem{thm}[lem]{Theorem}

\newtheorem{prop}[lem]{Proposition}

\newtheorem{cor}[lem]{Corollary}

\theoremstyle{definition}
\newtheorem{defn}[lem]{Definition}

\newtheorem{remark}[lem]{Remark}

\numberwithin{equation}{section}
\numberwithin{lem}{section}

\begin{document}

\title{$P$-adic incomplete gamma functions and Artin--Hasse-type series} 

\begin{abstract}
We define and study a $p$-adic analogue of the incomplete gamma function related to Morita's $p$-adic gamma function. We also discuss a combinatorial identity related to the Artin--Hasse series, which is a special case of the exponential principle in combinatorics. From this 
we deduce a curious $p$-adic property of $\#\mathrm{Hom}(G,S_n)$ for a topologically finitely generated group $G$, using a characterization of $p$-adic continuity for certain functions $f \colon \mathbb Z_{>0} \to \mathbb Q_p$ due to O'Desky--Richman. 
In the end, we give an exposition of some standard properties of the Artin--Hasse series. 
\end{abstract}


\author{Xiaojian Li}
\address{University of California, Berkeley}
\email{xiaojian-li@berkeley.edu}

\author{Jay Reiter}
\address{University of Illinois at Urbana-Champaign}
\email{jayr2@illinois.edu}

\author{Shiang Tang}
\address{Purdue University}
\email{shiangtang1989@gmail.com}

\author{Napoleon Wang}
\address{University of Arizona}
\email{napoleonw528@math.arizona.edu}

\author{Jin Yi}
\address{University of Illinois at Urbana-Champaign}
\email{hyunjin5@illinois.edu}

\maketitle 

\section{Introduction}
The first theme of this paper is the $p$-adic analogues of various gamma functions. 
Morita's $p$-adic gamma function \cite{morita} is a $p$-adic analogue of the classical gamma function that has been extensively studied in the literature. On the other hand, recent work of O'Desky--Richman \cite{odesky-richman} defines and studies a $p$-adic analogue of the \emph{incomplete gamma function}
\[
\Gamma(s,z) : = \int_{z}^{\infty} t^{s-1} e^{-t} dt
\]
defined for certain complex numbers $s,z$. However, it is unclear how their gamma function may be related to Morita's gamma function. In Section \ref{sec:incom gamma}, we define and study another $p$-adic analogue of the incomplete gamma function (Definition \ref{def:p adic incom gamma}, Theorem \ref{thm:p adic incom gamma}), which satisfies a recurrence relation similar to that of Morita's gamma function. 

The second theme of this paper is the Artin--Hasse series and related topics. The classical Artin--Hasse series 
\[
E_p(x) : = \exp \left ( \sum_{n \geq 0} \frac{x^{p^n}}{p^n} \right )
\]
has found numerous applications in number theory, $p$-adic analysis and combinatorics. An interesting combinatorial property is that the $n$-th EGF coefficient of $E_p(x)$ equals the number of elements in $S_n$ with $p$-power order. In Section \ref{sec:wohl}, we discuss a well-known generalization of this fact (Theorem \ref{thm:wohl id}), which is a special case of the exponential principle in combinatorics (see for example \cite{dress-thomas}). To keep this paper as self-contained as possible, we give a direct proof of this special case, which does not seem to be easily accessible in the literature. Combining Theorem \ref{thm:wohl id} with a characterization for $p$-adic continuity in terms of certain auxiliary quantities associated to a function $f \colon \mathbb Z_{>0} \to \mathbb Q_p$ \cite[Theorem 1.2]{odesky-richman}, we deduce that for a topologically finitely generated group $G$, the function $n \mapsto \#\mathrm{Hom}(G,S_n)$ (where the homomorphisms are assumed to be continuous w.r.t. the topology on $G$ and the discrete topology on $S_n$) can be extended to a continuous $p$-adic function on $\mathbb Z_p$ if and only if the number of open subgroups of $G$ with index $p$ is divisible by $p$ (Corollary \ref{cor:symm group p adic}). This generalizes an observation of O'Desky--Richman: the EGF coefficients of $E_p(x)$ can be extended to a continuous $\ell$-adic function on $\mathbb Z_{\ell}$ if and only if $p \neq \ell$. 

In Section \ref{sec:ah series}, we give an account of some standard facts on the Artin--Hasse series. For example, 
we give a proof of the fact that $E_p$, as an analytic function on $\mathfrak m_p$ (the unit open disk in $\mathbb C_p$) runs over the $p$-power roots of unity in $\mathbb C_p$. This is deduced from the fact that $E_p$ is a bijective isometry from $\mathfrak m_p$ to $1+\mathfrak m_p$ (Proposition \ref{prop:ah surj}), which seems to be well-known but we cannot find a proof in the literature.  

\textbf{Acknowledgements}: This paper is the end product of an Illinois Geometry Lab project (an undergraduate research project) conducted by the third author (S.T.) in Spring 2021. IGL research is supported by the Department of Mathematics at the University of Illinois at Urbana-Champaign. We would like to express deep gratitude to Ravi Donepuli, who was the Graduate Assistant of this project, for his constant support throughout the semester, without which the project would not have been carried out smoothly and successfully. S.T. would like to thank Andy O'Desky for helpful comments on an earlier version of this draft, and we thank the anonymous referee for comments and corrections that helped us improve the exposition.  

\subsection{Notation}
Let $p$ be an odd prime. 
We write $\mathbb Q_p$ for the field of $p$-adic numbers, write $\mathbb Z_p$ for the ring of $p$-adic integers, and write $\mathbb C_p$ for the completion of the algebraic closure of $\mathbb Q_p$ (``the field of $p$-adic complex numbers''). We write $|.|$ for the $p$-adic absolute value on $\mathbb C_p$, normalized so that $|p|=1/p$. For $a \in \mathbb C_p$ and $r>0$, let $B_a(r) : = \{ x \in \mathbb C_p \colon |x-a|<r \}$ be the open disk at $a$ of radius $r$. We write $\mathfrak m_p$ for $B_0(1) \subset \mathbb C_p$. Let $\exp_p \colon B_0((\frac{1}{p})^{1/(p-1)}) \xrightarrow{\sim} B_1((\frac{1}{p})^{1/(p-1)})$ be the $p$-adic exponential function with inverse function $\log_p$. 

\section{A $p$-adic incomplete gamma function} \label{sec:incom gamma}

\subsection{The classical incomplete gamma function}
Recall the classical incomplete gamma function:
\[ \Gamma(s,z) : = \int_z^{\infty} t^{s-1} e^{-t} dt \]
where $s$ is a complex number with positive real part, and $z$ is a real number. 
The following facts are well-known:
\begin{prop} \label{prop:classical incom gamma}
We have
\begin{itemize}
    \item $\lim_{z \to 0+} \Gamma(s,z) = \Gamma(s)$. 
    \item $\Gamma(1,z)=e^{-z}$.
    \item $\Gamma(s+1,z)=s \Gamma(s,z) + z^s e^{-z}$.
    \item For $s \in \mathbb Z_{>0}$,
    \[ \Gamma(s,z) = e^{-z} (s-1)! \sum_{k=0}^{s-1} \frac{z^k}{k!}. \]
\end{itemize}
\end{prop}

The incomplete function can be extended to a (multi-valued) holomorphic function in $s$ and $z$, but we do not discuss it here, see for example \cite[Section 5.1]{odesky-richman} for a discussion of this. 

\subsection{Morita's $p$-adic gamma function}
Morita's $p$-adic gamma function is the unique continuous function
$\Gamma_p(s)$ with $s \in \mathbb Z_p$ such that
\[ \Gamma_p(n) = (-1)^n \prod_{1 \leq k < n, p \nmid k} k \]
for $n \in \mathbb Z_{>0}$.
See \cite[Sections 35-39]{schikhof} for a detailed discussion of this function. 

For use of later discussions, we recall the following standard properties of $\Gamma_p$:

\begin{prop} \cite[Proposition 35.3]{schikhof} \label{prop:morita} 
Let $p$ be an odd prime. Then 
\begin{itemize}
    \item For all $s \in \mathbb Z_p$, $\Gamma_p(s+1)=h_p(s)\Gamma_p(s)$, where 
    $h_p(s)=-s$ if $s \in \mathbb Z_p^{\times}$, and $h_p(s)=-1$ if $s \in p\mathbb Z_p$.
    \item For all $x,y \in \mathbb Z_p$, $|\Gamma_p(x)-\Gamma_p(y)| \leq |x-y|$.
    \item $|\Gamma_p(s)|=1$ for all $s \in \mathbb Z_p$. 
\end{itemize}
\end{prop}

\subsection{A $p$-adic incomplete gamma function}
In this section, we define a $p$-adic analogue of the classical incomplete gamma function, which is closely related to Morita's $p$-adic gamma function. 

The following proposition is extracted from \cite[Theorem 34.1]{schikhof} and its proof. Let $K$ be a complete non-archimedean nontrivially valued field (e.g. $\mathbb Q_p$ or $\mathbb C_p$). 
\begin{prop} \label{prop:extension sum}
Let $f \colon \mathbb Z_p \to K$ be a continuous function. Then there is a unique continuous function $F \colon \mathbb Z_p \to K$ such that 
\begin{itemize}
    \item $F(s+1)=F(s)+f(s)$, $s \in \mathbb Z_p$.
    \item $F(0)=0$.
\end{itemize}
Moreover, let $||f|| : = \sup \{ |f(s)| \colon s \in \mathbb Z_p \}$, 
and let 
$\rho_i : = \sup \{ |f(s)-f(t)| \colon |s-t| \leq p^{-i} \}$ for $i \in \mathbb Z_{> 0}$. Then for integers $j \geq i \geq 1$ and $n \geq 1$, we have
\[
|F(n+p^j)-F(n)| \leq \max \{ \rho_i, p^{i-j} ||f|| \}. 
\]
\end{prop}

\begin{cor} \label{cor:extension sum}
Let $\{a_j\}$ be a sequence of elements in $K$. If $k \mapsto a_k$ extends to a continuous $p$-adic function on $\mathbb Z_p$, so does the partial sum $n \mapsto \sum_{k=0}^{n-1} a_k$.
Writing $\sum_{k=0}^{s-1} a_k$ with $s \in \mathbb Z_p$ for this extension, we have 
\[ \sum_{k=0}^{s-1} a_k = \lim_{j \to \infty} \sum_{k=0}^{n_j-1} a_k \]
where $\{ n_j \}$ is a sequence in $\mathbb Z_{>0}$ converging to $s \in \mathbb Z_p$ and the summation on the right hand side is the usual summation.
\end{cor}

\begin{proof}
Let $f \colon \mathbb Z_p \to K$ be the continuous extension of $a_k$ to $\mathbb Z_p$, and let $F$ be the associated function given by Proposition \ref{prop:extension sum}. The two bulleted points in Proposition \ref{prop:extension sum} imply that for $n \in \mathbb Z_{>0}$, $F(n)=\sum_{k=0}^{n-1} f(k)=\sum_{k=0}^{n-1} a_k$, which proves the first part of the corollary.
The second part is clear by continuity. 
\end{proof}

\begin{defn} \label{def:p adic incom gamma}
Define the \emph{$p$-adic incomplete gamma function} to be
\[
\Gamma_p(s,z) : = \exp_p(pz) \Gamma_p(s) \sum_{k=0}^{s-1} \frac{z^k}{\Gamma_p(k+1)} 
\]
where $s \in \mathbb Z_p$ and $z \in 1+p\mathbb Z_p$.
\end{defn}

Note the resemblance between this definition and the last formula in Proposition \ref{prop:classical incom gamma}: $k!=\Gamma(k+1)$ is replaced by $\Gamma_p(k+1)$ (similarly for $(s-1)!$), and the exponential function is replaced by the $p$-adic exponential function. 

\begin{lem}
Fix $z \in 1+p\mathbb Z_p$. The sequence $k \mapsto \frac{z^k}{\Gamma_p(k+1)}$ extends to a continuous $p$-adic function on $\mathbb Z_p$. It follows that $\Gamma_p(s,z)$ is a well-defined function on $\mathbb Z_p \times (1+p\mathbb Z_p)$. 
\end{lem}

\begin{proof}
By Proposition \ref{prop:morita}, $k \mapsto \Gamma_p(k+1)$ extends to a continuous function on $\mathbb Z_p$. On the other hand, \cite[Theorem 32.4]{schikhof} implies that for $z \in 1+p\mathbb Z_p$ (which is equivalent to $z$ being positive in the sense of loc. cit.) $z^k$ extends to a continuous function on $\mathbb Z_p$. By Corollary \ref{cor:extension sum} and Definition \ref{def:p adic incom gamma}, it follows that $\Gamma_p(s,z)$ is well-defined. 
\end{proof}

\begin{thm} \label{thm:p adic incom gamma}
The function $\Gamma_p(s,z)$ is the unique continuous function on $\mathbb Z_p \times (1+p\mathbb Z_p)$ satisfying 
\begin{enumerate}
    \item $\Gamma_p(1,z)=\exp_p(pz)$.
    \item $\Gamma_p(s+1,z)=h_p(s) \Gamma_p(s,z) + z^s \exp_p(pz)$ (where $h_p(s)$ is defined in Proposition \ref{prop:morita}). 
\end{enumerate}
\end{thm}

\begin{proof}
First we check that $\Gamma_p(s,z)$ satisfies (1)-(2). The first one is obvious. For the second one, we compute
\begin{align*}
    \Gamma_p(s+1,z) &= \exp_p(pz) \Gamma_p(s+1) \sum_{k=0}^{s} \frac{z^k}{\Gamma_p(k+1)} \\
    &=\exp_p(pz) h_p(s)\Gamma_p(s) \left( \sum_{k=0}^{s-1} \frac{z^k}{\Gamma_p(k+1)}+\frac{z^s}{\Gamma_p(s+1)} \right) \\
    &=h_p(s)\Gamma_p(s,z) + z^s \exp_p(pz)
\end{align*}
where the second equality follows from Propositions \ref{prop:morita} and \ref{prop:extension sum}.

Now we prove the continuity of $\Gamma_p(s,z)$. For this, it suffices to show that 
\[ F(s,z) : = \sum_{k=0}^{s-1} \frac{z^k}{\Gamma_p(k+1)} \]
is continuous on $\mathbb Z_p \times (1+p\mathbb Z_p)$. 

\emph{Claim}: if $\{s_n\}$ is a sequence of natural numbers converging to $s$, then $F(s_n,z)$ converges to $F(s,z)$ uniformly in $z$. 

Let $f_z(s) : = \frac{z^s}{\Gamma_p(s+1)} $. 
By Cauchy's criterion and the ultrametric inequality, it suffices to show that $|F(s_n,z)-F(s_{n-1},z)|$ approaches 0 uniformly in $z$ as $n \to \infty$. We may assume that $s_n=s_{n-1}+p^j$ for some integer $j>0$ (in particular $j$ depends on $n$, and goes to infinity with $n$).  
By the second part of Proposition \ref{prop:extension sum}, for any positive integer $i \leq j$, 
\begin{equation} \label{eq1}
  |F(s_n,z)-F(s_{n-1},z)| \leq \max \{ \rho_i, p^{i-j} ||f_z|| \}.  
\end{equation}
Note that $||f_z||=\sup \{ \frac{|z^s|}{|\Gamma_p(s+1)|} \colon s \in \mathbb Z_p \}$, which equals 1 by the assumption that $z \in 1+p\mathbb Z_p$ and the last point of Proposition \ref{prop:morita}. Also, 
\begin{align*}
    |z^s-z^t|&=|\exp_p(s\log_p(z))-\exp_p(t\log_p(z))| \\
    &=|s-t||\log_p(z)| \leq |s-t|
\end{align*}
where we have used the fact that $\exp_p \colon p\mathbb Z_p \to 1+p\mathbb Z_p$ is a bijective isometry with inverse $\log_p$. Combined with the fact that $1/\Gamma_p(s+1)$ is continuous, we see that $\rho_i$ (the constant associated to $f_z$ as in Proposition \ref{prop:extension sum}) approaches to 0 uniformly in $z$ as $i \to \infty$. It follows that (by taking $i=[j/2]$) the right hand side of (\ref{eq1}) approaches 0 as $n \to \infty$, which is uniform in $z$. This proves the claim.

From the claim it follows that 
\begin{itemize}
    \item For any $s \in \mathbb Z_p$, $\lim_{t \to s} F(t,z) = F(s,z)$ and the convergence is uniform in $z$. 
    \item For any $s \in \mathbb Z_p$, $z \mapsto F(s,z)$ is continuous (being the uniform limit of polynomial functions in $z$).
\end{itemize}
 
Finally, we show that $F$ is continuous at every $(s,z) \in \mathbb Z_p \times (1+p\mathbb Z_p)$. Fix $(s,z)$, let $(t,w)$ be a point in a neighborhood of $(s,z)$. We have
\[
|F(t,w)-F(s,z)| \leq \max \{ |F(t,w)-F(s,w)|, |F(s,w)-F(s,z)| \},
\]
$|F(t,w)-F(s,w)|$ (resp. $|F(s,w)-F(s,z)|$) can be made arbitrarily small by the first (resp. second) point above. 

Finally, we show that $\Gamma_p(s,z)$ is the unique continuous function satisfying (1)-(2). In fact, (1) and (2) determine $\Gamma_p(s,z)$ for $s \in \mathbb Z_{>0}$ and $z \in 1+p\mathbb Z_p$. Since $\mathbb Z_{>0}$ is dense in $\mathbb Z_p$, this determines $\Gamma_p(s,z)$ for all $s \in \mathbb Z_p$ by continuity. 
\end{proof}

\begin{remark}
Note that Theorem \ref{thm:p adic incom gamma}, (1) and (2) resemble the corresponding properties of the archimedean incomplete gamma function in Proposition \ref{prop:classical incom gamma}. It is natural to ask if there is an analogue of the identity $\lim_{z \to 0+} \Gamma(s,z) = \Gamma(s)$. Note that one can certainly extend $\Gamma_p(s,z)$ to a continuous function on $\mathbb Z_p \times \mathbb Z_p$ which specializes to $\Gamma_p(s)$ at $z=0$: for example, define $\Gamma_p(s,z)=0$ for $(s,z) \in \mathbb Z_p \times (k+p\mathbb Z_p)$ with $1<k<p$, and define $\Gamma_p(s,z)=\Gamma_p(s)$ for $(s,z) \in \mathbb Z_p \times p\mathbb Z_p$. But such an extension is unnatural at least because it is incompatible with the natural extension of $\Gamma_p(s,z)$ from $\mathbb Z_{>0} \times (1+p\mathbb Z_p)$ to $\mathbb Z_{>0} \times \mathbb Z_p$ given by the defining formula in Definition \ref{def:p adic incom gamma} (for $s \in \mathbb Z_{> 0}$, the summation $\sum_{k=0}^{s-1} \frac{z^k}{\Gamma_p(k+1)}$ makes sense for any $z$).
It is unclear to us whether there is a natural extension of $\Gamma_p(s,z)$ that recovers Morita's gamma function at $z=0$. 
\end{remark}

\textbf{Questions}: Is $\Gamma_p(s,z)$ analytic/locally analytic in $s$ (resp. $z$)? Does it have finitely many zeros on $\mathbb Z_p \times (1+p\mathbb Z_p)$?

\begin{remark}
O'Desky and Richman \cite[Theorem 1.1]{odesky-richman} defined another $p$-adic analog of the classical incomplete gamma function, which we denote by $\Gamma_p^*(s,z)$. It is a continuous function of $(s,z) \in \mathbb Z_p \times (1+p\mathbb Z_p)$ defined by the formula
\[
\Gamma_p^*(s,z) : = \exp_p(pz) \sum_{k \geq 0} z^{s-1-k} k! \binom{s-1}{k}
\]
and it satisfies 
\begin{itemize}
    \item $\Gamma_p^*(1,z)=\exp_p(pz)$.
    \item $\Gamma_p^*(s+1,z)=s \Gamma_p^*(s,z) + z^s \exp_p(pz)$.
\end{itemize}
Moreover, $\Gamma_p^*(s,z)$ is related to certain combinatorial objects, as explained in \emph{loc. cit.} 
\end{remark}


\section{A combinatorial identity} \label{sec:wohl}
Let $G$ be a topologically finitely generated group. For a positive integer $n$,
let $t_n$ be the number of continuous homomorphisms from $G$ to $S_n$ (the symmetric group on $n$ letters), and 
let $M_n$ be the number of open subgroups with index $n$. 
Note that $t_n$ and $M_n$ are necessarily finite.

The following theorem is a mild generalization of a theorem of Wohlfahrt \cite[Satz 1]{wohl} to topological groups: 

\begin{thm} \label{thm:wohl id}
\[
\exp \left ( \sum_{H \leq G} \frac{x^{[G:H]}}{[G:H]} \right )
= \sum_{n \geq 0} \frac{t_n}{n!} x^n
\]
where $H$ runs over open subgroups of $G$ with finite index.
\end{thm}

Identities of the above type have been extensively studied in group theory and combinatorics, see \cite[Proposition 1.1]{dress-thomas} for a general statement of this type. Since the proof of \emph{loc. cit.} falls out of an involved combinatorial theory, we have chosen to include a direct proof of the above theorem here, following the argument in \cite{wohl}. For a different proof of the special case when $G=\mathbb Z_p$, see \cite[Lemma 2.11]{kracht}.

\begin{lem}
Let $s_n$ be the number of continuous homomorphisms from $G$ to $S_n$ such that the image is a transitive subgroup of $S_n$.
Then $s_n=(n-1)!M_n$.
\end{lem}

\begin{proof}
For any open subgroup $H \subset G$ with index $n$, the natural action of $G$ on $G/H$ by left multiplication is transitive. 
Given a continuous, transitive action of $G$ on $\{1,2, \cdots, n\}$, the stablizer of 1 is an open subgroup $H$ of $G$, and this action may be identified with the canonical action of $G$ on $G/H$ by choosing a bijection between $\{1,2, \cdots, n\}$ and $G/H$ for which $1 \mapsto 1 \cdot H$. For each $H$, there are $(n-1)!$ such bijections. 
\end{proof}

\begin{lem}
Set $t_0=1$. We have
\[ t_n=\sum_{k=1}^{n} \binom{n-1}{k-1} s_k t_{n-k}. \]
\end{lem}
\begin{proof}
Given a continuous action of $G$ on $\{1,2, \cdots, n\}$, let $O(1)=\{1, i_2, i_3, \cdots , i_k\}$ be the orbit of 1. Then $G$ acts transitively on $O(1)$ and preserves $\{1,2, \cdots, n\} - O(1)$. Fix $k$, we can pick the $k-1$ elements $\{i_2, i_3, \cdots , i_k\}$ from $\{2, \cdots , n\}$ in $\binom{n-1}{k-1}$ ways, and for each choice, there are $s_k \cdot t_{n-k}$ possible actions of $G$ of the above type. The lemma follows by summing over $k=1,2, \cdots, n$. 
\end{proof}

The following proposition is essentially \cite[Satz 1]{wohl}: 
\begin{prop} \label{prop:wohl}
\[ t_n=(n-1)! \sum_{k=1}^n \frac{t_{n-k}M_k}{(n-k)!} = 
(n-1)! \left( t_0M_n+\frac{t_1M_{n-1}}{1!}+\cdots +\frac{t_{n-1}M_1}{(n-1)!} \right).
\]
\end{prop}

\begin{proof}
This follows immediately from the above two lemmas. 
\end{proof}

We will deduce Theorem \ref{thm:wohl id} from Proposition \ref{prop:wohl}. For this, we need the complete Bell polynomials, see \cite[p.205]{andrews}. 

\begin{defn}
The \emph{complete Bell polynomial} $B_n(x_1, \cdots, x_n)$
is defined by $B_0=1$, and 
\[
B_{n+1}(x_1, \cdots, x_{n+1})=\sum_{k=0}^n \binom{n}{k} B_{n-k}(x_1, \cdots, x_{n-k}) x_{k+1}.
\]
\end{defn}

For example, $B_1=x_1$, $B_2=x_1^2+x_2$, $B_3=x_1^3+3x_1x_2+x_3$.

The following combinatorial property is standard:

\begin{prop} \label{prop:bellpoly}
\[
\exp \left( \sum_{n \geq 1} x_n \frac{X^n}{n!} \right)
=\sum_{n \geq 0} B_n(x_1, \cdots, x_n) \frac{X^n}{n!}.
\]
\end{prop}

\begin{proof}[Proof of Theorem \ref{thm:wohl id}]
First we show that
\[
t_n=B_n(0!M_1,1!M_2,\cdots,(n-1)!M_n).
\]
We proceed by induction on $n$. The base case holds since $t_0=B_0=1$. Suppose the equality holds for $0,1,\cdots, n$. We have
\begin{align*}
    B_{n+1}(0!M_1,1!M_2,\cdots,n!M_{n+1}) 
    &=\sum_{k=0}^n \binom{n}{k} B_{n-k}(0!M_1,\cdots,(n-k-1)!M_{n-k}) \cdot (k!M_{k+1}) \\
    &=\sum_{k=0}^n \frac{n!}{(n-k)!}t_{n-k}M_{k+1} \\
    &=\sum_{k=1}^{n+1} \frac{n!}{(n-k+1)!}t_{n-k+1}M_{k}
    =t_{n+1}
\end{align*}
where the second equality follows from the inductive hypothesis, and the last equality follows from Proposition \ref{prop:wohl}.
Now
\begin{align*}
\exp \left ( \sum_{H \leq G} \frac{x^{[G:H]}}{[G:H]} \right )
&=\exp \left( \sum_{k \geq 1} \frac{M_k}{k} x^k \right) \\
&=\exp \left( \sum_{k \geq 1} M_k (k-1)! \frac{x^k}{k!}  \right) \\
&=\sum_{n \geq 0} B_n(0!M_1,1!M_2,\cdots,(n-1)!M_n) \frac{x^n}{n!} \\
&=\sum_{n \geq 0} \frac{t_n}{n!} x^n.
\end{align*}
where the third equality follows from Proposition \ref{prop:bellpoly}.
\end{proof}

Theorem \ref{thm:wohl id} implies the following curious property, which appears to be new: 
\begin{cor} \label{cor:symm group p adic}
The function $n \mapsto \# \mathrm{Hom} (G, S_n)$ extends to a continuous $p$-adic function on $\mathbb Z_p$ if and only if the number of open subgroups of $G$ with index $p$ is divisible by $p$. In particular, if $G$ is finite of order prime to $p$, then $n \mapsto \# \mathrm{Hom} (G, S_n)$ extends to a continuous $p$-adic function on $\mathbb Z_p$. 
\end{cor}

\begin{proof}
This relies on a characterization (due to O'Desky--Richman) for $p$-adic continuity in terms of certain auxiliary quantities associated to a function $f \colon \mathbb Z_{>0} \to \mathbb Q_p$. By \cite[Theorem 1.2]{odesky-richman} and Theorem \ref{thm:wohl id}, $t_n$ can be extended to a continuous $p$-adic function on $\mathbb Z_p$ if and only if $M_p \equiv M_1-1 \mod p$. Since $M_1=1$, this is equivalent to $p | M_p$. (Note that their term ``$p$-adically continuous'' is equivalent to ``extendable to a continuous $p$-adic function on $\mathbb Z_p$''.)
\end{proof}

The special case when $G=\mathbb Z_p$ is \cite[Corollary 1.3]{odesky-richman}.   

\section{The Artin--Hasse series} \label{sec:ah series}

The Artin--Hasse series is by definition the formal power series obtained through the composition of two power series
\[
E_p(x) : = \exp \left ( \sum_{n \geq 0} \frac{x^{p^n}}{p^n} \right )
\]
By Theorem \ref{thm:wohl id} (taking $G=\mathbb Z_p$), 
\[
E_p(x)=\sum_{n \geq 0} \frac{t_{p,n}}{n!} x^n
\]
where $t_{p,n}$ is the number of elements in $S_n$ with order a power of $p$. 

The following proposition is well-known, see the proof of \cite[Theorem 2.10]{kracht}.
\begin{prop} \label{prop:ah integral}
The coefficients $\frac{t_{p,n}}{n!}$ are $p$-integral. In particular, $E_p(x)$ converges for $x \in \mathfrak m_p$. 
\end{prop}


There are two other proofs of the above proposition in the literature: one uses an infinite product decomposition of $E_p(x)$, the other uses Dwork's lemma, see for example, \cite[Section 2]{conrad}. 

The following proposition is presumably well-known but we could not find an explicit statement in the literature:
\begin{prop} \label{prop:ah surj}
$x \mapsto E_p(x)$ defines a bijective isometry $\mathfrak m_p \to 1+\mathfrak m_p$. 
\end{prop}

\begin{proof}
The fact that $E_p$ is an isometry is standard, see \cite[Theorem 2.5]{conrad}. The new aspect in this proof is the argument for surjectivity. 

First we claim that for $0<r<1$,
\begin{equation} \label{eq2}
    \sup \left\{ \left | \frac{E_p(x)-E_p(y)}{x-y}-1 \right | \colon x \neq y, |x|, |y| \leq r \right\} \leq r.
\end{equation}
In fact,
\[
\left | \frac{E_p(x)-E_p(y)}{x-y}-1 \right |
    =\left | \frac{1}{x-y} \sum_{n \geq 2} \frac{t_{p,n}}{n!} (x^n-y^n) \right | 
    =\left | \sum_{n \geq 2} \frac{t_{p,n}}{n!} \left( \sum_{i=0}^{n-1} x^i y^{n-1-i} \right) \right | 
\]
By Proposition \ref{prop:ah integral}, $|\frac{t_{p,n}}{n!}| \leq 1$. Also, $|\sum_{i=0}^{n-1} x^i y^{n-1-i}|\leq r^{n-1}$ since $|x|,|y| \leq r$. Since $r<1$ and $n \geq 2$, it follows that the above quantity is at most $r$. 

(\ref{eq2}) and the fact that $|x+y|=\max(|x|,|y|)$ when $|x| \neq |y|$ imply that 
$\left | \frac{E_p(x)-E_p(y)}{x-y} \right |=1, \forall x,y \in \mathfrak m_p$ (by chosing appropriate $r$ for each pair of $x,y$), i.e.
$E_p$ is an isometry. 
In particular, $|E_p(x)-1|=|x|$ for all $x \in \mathfrak m_p$ and hence $E_p(x) \in 1+\mathfrak m_p$. It remains to show that $E_p$ maps \emph{onto} $ 1+\mathfrak m_p $. Let $c $ be a point in $1+\mathfrak m_p$ and choose $0<r<1$ such that $|c-1| \leq r$. 
Consider the function
\[
g(x) : = x-(E_p(x)-c), |x| \leq r.
\]
Observe that $|g(x)| \leq r$. In fact, $|E_p(x)-c| \leq \max \{ |E_p(x)-1|, |1-c|  \} = \max \{ |x|, |1-c|  \} \leq r$.
So $g$ is a continuous function from $\{|x| \leq r\}$ to itself. Moreover, 
\[
|g(x)-g(y)|=|x-y-(E_p(x)-E_p(y))|=|x-y| \left | \frac{E_p(x)-E_p(y)}{x-y}-1 \right | \leq r |x-y|
\]
where the last inequality follows from (\ref{eq2}).
So $g$ is a contraction on $\{|x| \leq r\}$, which is a complete metric space (by the completeness of $\mathbb C_p$). 
By Banach's contraction theorem (see for example, \cite[Appendix A.1]{schikhof}), $g$ has a fixed point in $\{|x| \leq r\}$, say $b$. Then $E_p(b)=c$, which proves the surjectivity of $E_p$. 
\end{proof}

\begin{cor} \label{cor:p power roots}
$E_p(x)$ runs over the $p$-power roots of unity in $\mathbb C_p$.
More precisely, for a primitive $p^k$-th root of unity $\zeta$, there exists $\alpha \in \mathbb C_p$ with $|\alpha|=(\frac{1}{p})^{1/(p-1)p^{k-1}}$ such that $E_p(\alpha)=\zeta$.
\end{cor}

\begin{proof}
Recall the following standard facts (see \cite[p.59, Lemma 10.1]{neukirch}):
$p=\prod_{1 \leq r <p^k, (r,p)=1} (1-\zeta^r)$
and $(1-\zeta^r)/(1-\zeta)$ is a unit for $(r,p)=1$. 
Taking $p$-adic absolute values on both sides, we find that
\[
\frac{1}{p} = |1-\zeta|^{\varphi(p^k)}= |1-\zeta|^{(p-1)p^{k-1}}
\]
and thus $|1-\zeta|=(\frac{1}{p})^{1/(p-1)p^{k-1}}$. In particular, $\zeta \in 1+\mathfrak m_p$, so Proposition \ref{prop:ah surj} implies that there exists $\alpha \in \mathbb C_p$ such that $E_p(\alpha)=\zeta$. We have $|\alpha|=|E_p(\alpha)-1|=|\zeta-1|=(\frac{1}{p})^{1/(p-1)p^{k-1}}$ (the first equality holds since $E_p$ is an isometry).
\end{proof}

For a different proof of this fact, see \cite{conrad}.

\begin{remark}
In contrast, the $p$-adic exponential function $\exp_p$ does not run over any $p$-power roots of unity in $\mathbb C_p$. In fact, the range of $\exp_p$ equals $B_1((\frac{1}{p})^{\frac{1}{p-1}})$ (see for example, \cite[Theorem 25.6 and Proposition 44.1]{schikhof}), which does not contain any nontrivial $p$-power roots of unity by the proof above. 
\end{remark}

\end{document}